\documentclass[11pt,reqno]{article}
\usepackage{amssymb,amscd,amsmath,amsthm}
\newtheorem{thm}{Theorem}[section]
\newtheorem{prop}[thm]{Proposition}
\newtheorem{lemma}[thm]{Lemma}
\newtheorem{cor}[thm]{Corollary}
\newtheorem{definition}[thm]{Definition}
\newtheorem{remark}[thm]{Remark}
\newtheorem{example}[thm]{Example}

\def\bC{\mathbb{C}}
\def\bR{\mathbb{R}}
\def\bN{\mathbb{N}}
\def\pos{\mathrm{pos}}
\def\mono{\mathrm{mono}}
\def\conv{\mathrm{conv}}
\def\eps{\varepsilon}
\def\<{\langle}
\def\>{\rangle}
\def\sign{\mathrm{sign}\,}
\def\Tr{\mathrm{Tr}\,}

\begin{document}

\ \vskip 1cm 
\centerline{\LARGE {\bf Monotonicity for entrywise functions of matrices}} 
\bigskip
\bigskip

\centerline{\Large
Fumio Hiai\footnote{Partially supported by Grant-in-Aid for Scientific Research
(B)17340043.}}

\medskip
\begin{center}
Graduate School of Information Sciences, Tohoku University \\
Aoba-ku, Sendai 980-8579, Japan \\
E-mail: hiai@math.is.tohoku.ac.jp
\end{center}

\bigskip
\begin{abstract}
We characterize real functions $f$ on an interval $(-\alpha,\alpha)$ for which the
entrywise matrix function $[a_{ij}]\mapsto[f(a_{ij})]$ is positive, monotone and
convex, respectively, in the positive semidefiniteness order. Fractional power
functions are exemplified and related weak majorizations are shown.


\end{abstract}

\bigskip
\section*{Introduction}

There are two important notions of order for matrices; one is the order induced by
positive semidefiniteness and the other is that induced by the positive cone of
entrywise nonnegative matrices. On the other hand, there are two ways in applying
functions (defined on an interval) to matrices, the usual functional calculus
$A\mapsto f(A)$ and the entrywise calculus $A\mapsto f[A]$. In this way, one may take
the following four combinations to study monotonicity or convexity for matrix
functions:
\begin{itemize}
\item functional calculus and positive semidefiniteness,
\item functional calculus and entrywise positivity,
\item entrywise calculus and positive semidefiniteness,
\item entrywise calculus and entrywise positivity.
\end{itemize}
The last situation is trivial; it has nothing to do with matrices. The first situation
is most standard and most important in matrix theory. We have a well-developed
theory of operator monotone and operator convex functions initiated by L\"owner; a
comprehensive exposition on the subject is found in \cite{Bh}. The second one was
treated by Hansen \cite{Ha}, and the third one is the subject of the present paper.

In \cite{FH} FitzGerald and Horn considered entrywise fractional power (i.e.,
fractional Hadamard or Schur power) $A^{(p)}:=[a_{ij}^p]$ for numbers $p>0$ and for
positive semidefinite and entrywise nonnegative matrices $A=[a_{ij}]$. They
characterized the numbers $p$ for which $A^{(p)}\ge0$ (positive semidefinite) for all
entrywise nonnegative $A\ge0$ and those $p$ for which $A^{(p)}\ge B^{(p)}$ for all
entrywise nonnegative $A\ge B\ge0$. These are typical results in the third situation
mentioned above, motivating us to consider the same problem in more general settings.

In this paper we treat a real function on an open interval $(-\alpha,\alpha)$ with
$0<\alpha\le\infty$. For a Hermitian complex matrix $A$ whose eigenvalues are in
$(-\alpha,\alpha)$ let $f(A)$ denote the usual functional calculus of $A$ by $f$. On
the other hand, for a real matrix $A=[a_{ij}]$ whose entries are all in
$(-\alpha,\alpha)$ we write $f[A]$ for the matrix obtained by applying $f$ entrywise
to $A$, i.e., $f[A]=[f(a_{ij})]$. Let $M_n(\bR)$ denote the set of $n\times n$ real
matrices. We say that $f$ is S-positive if $f[A]\ge0$ (positive semidefinite) for
every $A\ge0$ in $M_n(\bR)$ of any $n$ with entries in $(-\alpha,\alpha)$, and that
$f$ is S-monotone if $f[A]\ge f[B]$ (in the order of positive semidefiniteness) for
every $A\ge B\ge0$ in $M_n(\bR)$ of any $n$ with all entries in $(-\alpha,\alpha)$.
Moreover, $f$ is said to be S-convex if $A\mapsto f[A]$ satisfies the convexity
property for every pair $A\ge B\ge0$ as above. The main aim of this paper is to
completely characterize these three classes of functions on the interval
$(-\alpha,\alpha)$. But we also discuss the three classes for each fixed order $n$.

In \cite{Ha} a real function $f$ on $(-\alpha,\alpha)$ was said to be m-positive,
m-monotone and m-convex if it satisfies the properties similar to, respectively,
those of S-positive, S-monotone and S-convex functions in the second situation
mentioned above, i.e., in the setting of the usual functional calculus $f(A)$ and
the order of entrywise positivity. Rather surprisingly, Hansen's characterization
in \cite{Ha} is completely the same as ours in Theorem \ref{T-4.1}; thus the classes
of m-positive, m-monotone and m-convex functions on $(-\alpha,\alpha)$ coincide with,
respectively, those of S-positive, S-monotone and S-convex functions on
$(-\alpha,\alpha)$. Here it should be remarked that the whole structure of our proof
of Theorem \ref{T-4.1} is somewhat similar to that in \cite{Ha} though there are many
differences between the details of the two proofs.

The paper is organized as follows. The precise definitions of S-positive, S-monotone
and S-convex functions together with those for each fixed order $n$ are presented in
Section 1. In Section 2 we then obtain complete characterizations of the three classes
of the first non-trivial order $n=2$ as well as some necessary conditions for those of
the next order $n=3$. These discussions in lower order cases are indispensable in
proving our main theorem. In Section 3 we demonstrate several relations among the
three classes of order $n$ when $n$ varies. For instance, we show that S-positive of
order $2n$ implies S-monotone of order $n$ and S-monotone of order $2n$ implies
S-convex of order $n$. With some preparations in Sections 2 and 3 the main theorem
(Theorem \ref{T-4.1}) is proven in Section 4. Next in Section 5 we deal with examples
of fractional power functions and slightly extend the results in \cite{FH} mentioned
above. Finally in Section 6 we obtain related weak majorizations involving entrywise
matrix functions.

\section{Definitions}
\setcounter{equation}{0}

The set of complex $n\times n$ matrices is denoted by $M_n(\bC)$, and that of real
$n\times n$ matrices is by $M_n(\bR)$. The symbol $J$ stands for the $n\times n$
matrix with all entries equal to $1$ (i.e., the identity matrix for the Schur
product) while $I$ is the usual $n\times n$ identity matrix. For $A\in M_n(\bC)$,
$A\ge0$ means that $A$ is positive semidefinite. For Hermitian $A,B\in M_n(\bC)$
(in particular, for symmetric $A,B\in M_n(\bR)$), $A\ge B$ means $A-B\ge0$. Throughout
the paper we fix any $\alpha$ with $0<\alpha\le\infty$. For a real function $f$ on the
open interval $(-\alpha,\alpha)$ and for a matrix $A=[a_{ij}]$ in $M_n(\bR)$ such that
$a_{ij}\in(-\alpha,\alpha)$ for all $1\le i,j\le n$, we write $f[A]$ for the matrix
obtained by applying $f$ to all entries of $A$, i.e.,
$$
f[A]:=[f(a_{ij})].
$$

\begin{definition}{\rm
For a real function $f$ on $(-\alpha,\alpha)$ and for $n\in\bN$, we introduce the
following three notions:
\begin{itemize}
\item[\rm(i)] $f$ is {\it S-positive} (or {\it Schur positive}\,) {\it of order $n$}
if
$$
f[A]\ge0
$$
for all $A\ge0$ in $M_n(\bR)$ with entries in $(-\alpha,\alpha)$.
\item[\rm(ii)] $f$ is {\it S-monotone} (or {\it Schur monotone}\,) {\it of order $n$}
if
$$
A\ge B\ge0\ \ \Longrightarrow\ \ f[A]\ge f[B]
$$
for all symmetric $A,B\in M_n(\bR)$ with entries in $(-\alpha,\alpha)$.
\item[\rm(iii)] $f$ is {\it S-convex} (or {\it Schur convex}\,) {\it of order $n$}
if
$$
A\ge B\ge0\ \ \Longrightarrow
\ \ f[\lambda A+(1-\lambda)B]\le\lambda f[A]+(1-\lambda)f[B],\quad0\le\lambda\le1
$$
for all symmetric $A,B\in M_n(\bR)$ with entries in $(-\alpha,\alpha)$.
\end{itemize}
We denote by $S_\pos^{(n)}(-\alpha,\alpha)$, $S_\mono^{(n)}(-\alpha,\alpha)$ and
$S_\conv^{(n)}(-\alpha,\alpha)$ the classes of all real functions on
$(-\alpha,\alpha)$ which are respectively S-positive, S-monotone and S-convex of
order $n$. Moreover, when $f$ is S-positive (resp., S-monotone, S-convex) of all
order $n$, we say that $f$ is {\it S-positive} (resp., {\it S-monotone},
{\it S-convex}).
}\end{definition}

It is obvious that each class of $S_\pos^{(n)}(-\alpha,\alpha)$,
$S_\mono^{(n)}(-\alpha,\alpha)$ and $S_\conv^{(n)}(-\alpha,\alpha)$ becomes smaller
as $n$ increases. The classes $S_\pos^{(1)}(-\alpha,\alpha)$,
$S_\mono^{(1)}(-\alpha,\alpha)$ and $S_\conv^{(1)}(-\alpha,\alpha)$ are the sets of
real functions on $(-\alpha,\alpha)$ which are nonnegative, non-decreasing and
convex, respectively, on $[0,\alpha)$ in usual sense as real functions with no
requirement on $f|_{(-\alpha,0)}$. Also it is clear that if $f$ is S-monotone of
order $n$ and $f(0)\ge0$, then $f$ is S-positive of order $n$.

One might consider the condition $A\ge B\ge0$ in the above definitions (ii) and (iii)
too restrictive when these definitions are compared with those of operator monotonicity
and operator convexity. However, the next proposition shows that this restriction is
necessary when we treat entrywise matrix functions $[a_{ij}]\mapsto[f(a_{ij})]$.

\begin{prop}
Let $f$ be a real function on $(-\alpha,\alpha)$.
\begin{itemize}
\item[\rm(1)] If $A\ge B$ implies $f[A]\ge f[B]$ for all symmetric $A,B\in M_2(\bR)$
with entries in $(-\alpha,\alpha)$, then $f$ is affine on $(-\alpha,\alpha)$.
\item[\rm(2)] If $f[\lambda A+(1-\lambda)B]\le\lambda f[A]+(1-\lambda)f[B]$ for all
$0\le\lambda\le1$ and all symmetric $A,B\ge0$ in $M_2(\bR)$ with entries in
$(-\alpha,\alpha)$, then $f$ is affine on $(-\alpha,\alpha)$.
\end{itemize}
\end{prop}

\begin{proof}
(1)\enspace We may assume $f(0)=0$ by taking $f-f(0)$ instead of $f$. The assumption
of (1) obviously implies that $f$ is non-decreasing on $(-\alpha,\alpha)$; so
$f(x)\ge0$ for $0\le x<\alpha$ and $f(x)\le0$ for $-\alpha<x\le0$. Let $0\le a<\alpha$
and $0<\lambda<1$. Since
$$
\bmatrix a&\lambda a\\\lambda a&a\endbmatrix
\ge\bmatrix(1-\lambda)a&0\\0&(1-\lambda)a\endbmatrix,\quad
\bmatrix\lambda a&(1-\lambda)a\\(1-\lambda)a&\lambda
a\endbmatrix\ge\bmatrix0&a\\a&0\endbmatrix,
$$
we get
$$
\bmatrix f(a)&f(\lambda a)\\f(\lambda a)&f(a)\endbmatrix
\ge\bmatrix f((1-\lambda)a)&0\\0&f((1-\lambda)a)\endbmatrix,
$$
$$
\bmatrix f(\lambda a)&f((1-\lambda)a)\\
f((1-\lambda)a)&f(\lambda a)\endbmatrix
\ge\bmatrix0&f(a)\\f(a)&0\endbmatrix.
$$
From these we obtain $f(a)=f(\lambda a)+f((1-\lambda)a)$, which means that $f$ is
affine on $[0,\alpha)$. Furthermore, since $\bmatrix a&-a\\-a&a\endbmatrix\ge0$ and
$\bmatrix-a&a\\a&-a\endbmatrix\le0$, we have
$\bmatrix f(a)&f(-a)\\f(-a)&f(a)\endbmatrix\ge0$ and
$\bmatrix-f(-a)&-f(a)\\-f(a)&-f(-a)\endbmatrix\ge0$. These imply that $f(-a)=-f(a)$
for all $a\in[0,\alpha)$. So $f$ is affine on $(-\alpha,\alpha)$.

(2)\enspace Let $0<a<\alpha$ and $s,t\in[-a,a]$. Since
$\bmatrix a&s\\s&a\endbmatrix,\bmatrix a&t\\t&a\endbmatrix\ge0$, the assumption of (2)
implies that for every $0<\lambda<1$
$$
\bmatrix f(a)&f(\lambda s+(1-\lambda)t)\\
f(\lambda s+(1-\lambda)t)&f(a)\endbmatrix
\le\bmatrix f(a)&\lambda f(s)+(1-\lambda)f(t)\\
\lambda f(s)+(1-\lambda)f(t)&f(a)\endbmatrix
$$
and so we obtain $f(\lambda s+(1-\lambda)t)=\lambda f(s)+(1-\lambda)f(t)$. Hence $f$
is affine on $(-\alpha,\alpha)$.
\end{proof}

\begin{example}\label{E-1.3}{\rm
For each $k\in\bN$ consider the function $f(x)=x^k$ on $\bR$ and write $A^{(k)}:=f[A]$
for this $f$, that is, $A^{(k)}$ stands for the Schur product $A\circ\cdots\circ A$
($k$-fold). If $A\ge B\ge0$ in $M_n(\bR)$, then the well-known Schur theorem gives
$A^{(k)}\ge B^{(k)}$. So $f(x)=x^k$ is S-monotone (hence S-positive). This is the
reason why we use the term ``Schur positive," etc. Furthermore, $f(x)=x^k$ is
S-convex. This is trivial when $k=1$. Assuming the S-convexity of $x^k$, for
$A\ge B\ge0$ and $0\le\lambda\le1$ we have
\begin{align*}
(\lambda A+(1-\lambda)B)^{(k+1)}
&\le(\lambda A+(1-\lambda)B)\circ(\lambda A^{(k)}+(1-\lambda)B^{(k)}) \\
&=\lambda A^{(k+1)}+(1-\lambda)B^{(k+1)}
-\lambda(1-\lambda)(A-B)\circ(A^{(k)}-B^{(k)}) \\
&\le\lambda A^{(k+1)}+(1-\lambda)B^{(k+1)}
\end{align*}
by repeated use of the Schur theorem. Hence we get the S-convexity of $x^{k+1}$ as
well. Consequently, when $f$ has a series expansion
$f(x)=\sum_{k=0}^\infty\alpha_kx^k$ with radius $r>0$ of convergence, the function
$f$ on $(-r,r)$ is
\begin{itemize}
\item[(i)] S-positive if $\alpha_k\ge0$ for all $k\ge0$,
\item[(ii)] S-monotone if $\alpha_k\ge0$ for all $k\ge1$,
\item[(iii)] S-convex if $\alpha_k\ge0$ for all $k\ge2$.
\end{itemize}
}\end{example}

The main result shown in the present paper is that the sufficient conditions in the
above (i)--(iii) are even necessary, that is, the functions given above actually
exhaust the S-positive, S-monotone and S-convex functions, respectively.

\section{Lower order cases}
\setcounter{equation}{0}

The aim of this section is to give concrete descriptions of functions in the classes
$S_\pos^{(n)}(-\alpha,\alpha)$, $S_\mono^{(n)}(-\alpha,\alpha)$ and
$S_\conv^{(n)}(-\alpha,\alpha)$ in the lower order cases $n=2$ and $n=3$. This is the
first task to be done toward the goal of our characterization problem. 

Let $f$ be a nonnegative real function $f$ on the open interval $(0,\alpha)$. We say
that $f$ is {\it $\sqrt{\phantom x}$-submultiplicative} if
$$
f(\sqrt{st})\le\sqrt{f(s)f(t)}
\quad\mbox{for all }s,t\in(0,\alpha).
$$
The class of non-decreasing and $\sqrt{\phantom x}$-submultiplicative functions on
$(0,\alpha)$ is described as follows.

\begin{lemma}\label{L-2.1}
For a nonnegative function $f$ on $(0,\alpha)$ the following conditions are
equivalent:
\begin{itemize}
\item[\rm(a)] $f$ is non-decreasing and $\sqrt{\phantom x}$-submultiplicative;
\item[\rm(b)] $f$ is non-decreasing, continuous and
$\sqrt{\phantom x}$-submultiplicative;
\item[\rm(c)] $f$ is identically zero, or else there is a non-decreasing convex
function $g$ on $(-\infty,\log\alpha)$ such that $f(t)=\exp g(\log t)$ for all
$t\in(0,\alpha)$.
\end{itemize}
\end{lemma}

\begin{proof}
It is straightforward to see that (c) $\Rightarrow$ (b) $\Rightarrow$ (a). To prove
(a) $\Rightarrow$ (c), let $f$ be a non-decreasing and
$\sqrt{\phantom x}$-submultiplicative function on $(0,\alpha)$ which is not
identically zero. It is easily seen that $f(t)>0$ for all $t\in(0,\alpha)$. For every
$t\in(0,\alpha)$ and $0<\eps\le t/5$, since $t+\eps\le\sqrt{(t+4\eps)(t-\eps)}$, we get
$$
f(t+\eps)\le f(\sqrt{(t+4\eps)(t-\eps)})\le\sqrt{f(t+4\eps)f(t-\eps)}.
$$
Letting $\eps\searrow0$ gives
$$
\lim_{s\to t+0}f(s)\le\lim_{s\to t-0}f(s),
$$
which implies the continuity of $f$ at $t$. Now define a function $g$ on
$(-\infty,\log\alpha)$ by $g(x)=\log f(e^x)$ for $-\infty<x<\log\alpha$ so that
$f(t)=\exp g(\log t)$ for $0<t<\alpha$. Then $g$ is non-decreasing and continuous on
$(-\infty,\log\alpha)$ as so is $f$ on $(0,\alpha)$. The
$\sqrt{\phantom x}$-submultiplicativity and the continuity of $f$ imply the convexity
of $g$, and hence (c) follows.
\end{proof}

We denote by $\Phi(0,\alpha)$ the set of all nonnegative functions on $(0,\alpha)$
satisfying the equivalent conditions (a)--(c) in Lemma \ref{L-2.1}.

\begin{prop}\label{P-2.2}
For a real function $f$ on $(-\alpha,\alpha)$, $f\in S_\pos^{(2)}(-\alpha,\alpha)$ if
and only if $f|_{(0,\alpha)}\in\Phi(0,\alpha)$, $0\le f(0)\le f(0+)$
$(:=\lim_{t\searrow0}f(t))$ and $|f(-t)|\le f(t)$ for all $0<t<\alpha$.
\end{prop}

\begin{proof}
Assume $f\in S_\pos^{(2)}(-\alpha,\alpha)$. If $0\le t<s<\alpha$, then
$\bmatrix s&t\\t&s\endbmatrix\ge0$ implies
$\bmatrix f(s)&f(t)\\f(t)&f(s)\endbmatrix\ge0$, so $0\le f(t)\le f(s)$. Hence $f$ is
nonnegative and non-decreasing on $[0,\alpha)$. For all $s,t\in(0,\alpha)$, since
$\bmatrix s&\sqrt{st}\\\sqrt{st}&t\endbmatrix\ge0$, we get
$\bmatrix f(s)&f(\sqrt{st})\\f(\sqrt{st})&f(t)\endbmatrix\ge0$ and so
$f(\sqrt{st})\le\sqrt{f(s)f(t)}$, i.e., $f$ is $\sqrt{\phantom x}$-submultiplicative
on $(0,\alpha)$. Moreover, for $0<t<\alpha$, we get
$\bmatrix f(t)&f(-t)\\f(-t)&f(t)\endbmatrix\ge0$ so that $|f(-t)|\le f(t)$.

Conversely assume that $f$ satisfies the conditions stated in the proposition. Let
$\bmatrix a&c\\c&b\endbmatrix\ge0$ in $M_2(\bR)$ with $a,b,c\in(-\alpha,\alpha)$;
then $a,b\ge0$ and $c^2\le ab$. If $c=0$, then
$\bmatrix f(a)&f(0)\\f(0)&f(b)\endbmatrix\ge0$ since $0\le f(0)\le f(a),f(b)$. If
$c\ne0$, then $a,b>0$ and $f(c)^2\le f(|c|)^2\le f(\sqrt{ab})^2\le f(a)f(b)$, so
$\bmatrix f(a)&f(c)\\f(c)&f(b)\endbmatrix\ge0$. Hence
$f\in S_\pos^{(2)}(-\alpha,\alpha)$.
\end{proof}

We denote by $\Psi^{(1)}(-\alpha,\alpha)$ the set of all measurable real functions
$f$ on $(-\alpha,\alpha)$ such that $f|_{(0,\alpha)}\in\Phi(0,\alpha)$ and
$|f(-t)|\le f(t)$ for a.e.\ $t\in(0,\alpha)$ (with respect to the Lebesgue measure).
Moreover, let $\Psi^{(2)}(-\alpha,\alpha)$ denote the set of all continuous
functions $f$ on $(-\alpha,\alpha)$ which is differentiable on $(0,\alpha)$ and
differentiable a.e.\ on $(-\alpha,0)$ with $f'\in\Psi^{(1)}(-\alpha,\alpha)$.
In other words, $f\in\Psi^{(2)}(-\alpha,\alpha)$ if and only if there exists
$g\in\Psi^{(1)}(-\alpha,\alpha)$ such that
$$
f(t)-f(0)=\int_0^tg(s)\,ds\quad\mbox{for }-\alpha<t<\alpha.
$$

\begin{prop}\label{P-2.3}
$S_\mono^{(2)}(-\alpha,\alpha)=\Psi^{(2)}(-\alpha,\alpha)$.
\end{prop}

\begin{proof}
Assume $f\in S_\mono^{(2)}(-\alpha,\alpha)$. We may and do assume $f(0)=0$ by taking
$f-f(0)$ instead of $f$. First note that $f\in S_\pos^{(2)}(-\alpha,\alpha)$ and
hence $f$ is continuous on $(0,\alpha)$ by Proposition \ref{P-2.2}. For every
$0<t<s<\alpha$, since
$\bmatrix s&t\\t&s\endbmatrix\ge\bmatrix{s+t\over2}&{s+t\over2}\\
{s+t\over2}&{s+t\over2}\endbmatrix\ge0$, we get
$\bmatrix f(s)&f(t)\\f(t)&f(s)\endbmatrix\ge\bmatrix
f({s+t\over2})&f({s+t\over2})\\f({s+t\over2})&f({s+t\over2})\endbmatrix$.
By multiplying $\bmatrix1&1\\0&0\endbmatrix$ from the left and
$\bmatrix1&0\\1&0\endbmatrix$ from the right this implies
$f({s+t\over2})\le{f(s)+f(t)\over2}$. Hence $f$ is convex on $(0,\alpha)$; so it is
right-differentiable on $(0,\alpha)$ so that the right-derivative $f'_+$ is
non-decreasing on $(0,\alpha)$. For each $a,b,c\in(0,\alpha)$ with $c^2\le ab$
and for $\eps>0$ small enough, since $\bmatrix a+\eps&c+\eps\\c+\eps&b+\eps\endbmatrix
\ge\bmatrix a&c\\c&b\endbmatrix\ge0$, we get
$\bmatrix{f(a+\eps)-f(a)\over\eps}&{f(c+\eps)-f(c)\over\eps}\\
{f(c+\eps)-f(c)\over\eps}&{f(b+\eps)-f(b)\over\eps}\endbmatrix\ge0$.
Letting $\eps\searrow0$ gives
$\bmatrix f'_+(a)&f'_+(c)\\f'_+(c)&f'_+(b)\endbmatrix\ge0$ so that
$f'_+(c)^2\le f'_+(a)f'_+(b)$, i.e., $f'_+$ is $\sqrt{\phantom x}$-submultiplicative
on $(0,\alpha)$. Hence Lemma \ref{L-2.1} implies that $f'_+$ is continuous on
$(0,\alpha)$ and so $f$ is differentiable on $(0,\alpha)$ with
$f'|_{(0,\alpha)}\in\Phi(0,\alpha)$. Let $a\in(0,\alpha)$ and $0<\eps<\alpha-a$. Since
$\bmatrix a+\eps&a\\a&a\endbmatrix\ge\bmatrix\eps&0\\0&0\endbmatrix\ge0$,
we have $\bmatrix f(a+\eps)&f(a)\\f(a)&f(a)\endbmatrix\ge\bmatrix
f(\eps)&0\\0&0\endbmatrix$. Multiply $\bmatrix1&-1\\0&0\endbmatrix$ from the left and
$\bmatrix1&0\\-1&0\endbmatrix$ from the right to get $f(a+\eps)-f(a)\ge f(\eps)\ge0$.
Hence $\lim_{\eps\searrow0}f(\eps)=0=f(0)$ thanks to the continuity at $a$; so $f$ is
right-continuous at $0$ (the left-continuity at $0$ follows as well from the proof of
the absolute continuity of $f|_{(-\alpha,0]}$ below). Put $\tilde f(t):=f(-t)$ for
$0\le t<\alpha$. For every $0<a<\alpha$ and $0=t_0<t_1<\dots<t_k=a$, since
$\bmatrix t_j&-t_j\\-t_j&t_j\endbmatrix\ge
\bmatrix t_{j-1}&-t_{j-1}\\-t_{j-1}&t_{j-1}\endbmatrix\ge0$, we have
$\bmatrix f(t_j)&\tilde f(t_j)\\\tilde f(t_j)&f(t_j)\endbmatrix\ge
\bmatrix f(t_{j-1})&\tilde f(t_{j-1})\\\tilde f(t_{j-1})&f(t_{j-1})
\endbmatrix$. Therefore,
$$
|\tilde f(t_j)-\tilde f(t_{j-1})|\le f(t_j)-f(t_{j-1})\,,
\qquad 1\le j\le k,
$$
which yields
$$
\sum_{j=1}^k|\tilde f(t_j)-\tilde f(t_{j-1})|\le f(a).
$$
Hence $\tilde f$ is absolutely continuous on $[0,a]$ for any $0<a<\alpha$ and its total
variation on $[0,a]$ is dominated by $f(a)$ ($=\int_0^af'(t)\,dt$). This shows that
$\tilde f$ is differentiable a.e.\ on $(0,\alpha)$ and $|\tilde f'(t)|\le f'(t)$ for
a.e.\ $t\in(0,\alpha)$, that is, $f$ is differentiable a.e.\ on $(-\alpha,0)$ and
$|f'(-t)|\le f'(t)$ for a.e.\ $t\in(0,\alpha)$. Hence $f\in\Psi^{(2)}(-\alpha,\alpha)$.

Conversely assume $f\in\Psi^{(2)}(-\alpha,\alpha)$. Then there exists
$g\in\Psi^{(1)}(-\alpha,\alpha)$ such that
$$
f(t)-f(0)=\int_0^tg(s)\,ds\quad\mbox{for }-\alpha<t<\alpha.
$$
For any $s,t\in(-\alpha,\alpha)$ we have
$$
f(s)-f(t)=(s-t)\int_0^1g(\lambda s+(1-\lambda)t)\,d\lambda.
$$
Let $A\ge B\ge0$ in $M_2(\bR)$ with entries in $(-\alpha,\alpha)$. To prove
$f[A]\ge f[B]$, we may assume by continuity that $a_{ij}\ne b_{ij}$ for all $i,j$,
where $A=[a_{ij}]$ and $B=[b_{ij}]$. Under this assumption we have
$g[\lambda A+(1-\lambda)B]\ge0$ for a.e.\ $\lambda\in(0,1)$ as in the proof of
Proposition \ref{P-2.2}. So $f[A]\ge f[B]$ is obtained from the expression
$$
f[A]-f[B]=(A-B)\circ\int_0^1g[\lambda A+(1-\lambda)B]\,d\lambda.
$$
Hence $f\in S_\mono^{(2)}(-\alpha,\alpha)$.
\end{proof}

\begin{prop}\label{P-2.4}
For a real function $f$ on $(-\alpha,\alpha)$,
$f\in S_\conv^{(2)}(-\alpha,\alpha)$ if and only if $f$ is differentiable on
$(-\alpha,\alpha)$ and $f'\in\Psi^{(2)}(-\alpha,\alpha)$
$(=S_\mono^{(2)}(-\alpha,\alpha))$. Hence, if $f\in S_\conv^{(2)}(-\alpha,\alpha)$,
then $f$ is continuously differentiable on $(-\alpha,\alpha)$.
\end{prop}

\begin{proof}
By Proposition \ref{P-2.3} it suffices to prove the first assertion. Assume
$f\in S_\conv^{(2)}(-\alpha,\alpha)$. Obviously $f$ is convex on $[0,\alpha)$ and so
right-differentiable on $(0,\alpha)$. For $0\le t\le s<\alpha$ and $0<\eps<\alpha-s$,
since $\bmatrix s+\eps&s\\s&s\endbmatrix\ge\bmatrix t+\eps&t\\t&t\endbmatrix\ge0$, we
have for $0\le\lambda\le1$
\begin{align*}
&\bmatrix f(\lambda(s+\eps)+(1-\lambda)(t+\eps))
&f(\lambda s+(1-\lambda)t)\\f(\lambda s+(1-\lambda)t)
&f(\lambda s+(1-\lambda)t)\endbmatrix \\
&\qquad\le\lambda\bmatrix f(s+\eps)&f(s)\\f(s)&f(s)\endbmatrix
+(1-\lambda)\bmatrix f(t+\eps)&f(t)\\f(t)&f(t)\endbmatrix,
\end{align*}
which implies that
$$
f(\lambda s+(1-\lambda)t+\eps)-f(\lambda s+(1-\lambda)t) \\
\le\lambda(f(s+\eps)-f(s))+(1-\lambda)(f(t+\eps)-f(t)).
$$
By dividing by $\eps$ and then letting $\eps\searrow0$ we see that $f'_+$ is convex
on $(0,\alpha)$ and so it is continuous on $(0,\alpha)$. Now let
$f_0(t):={f(t)+f(-t)\over2}$ and $f_1(t):={f(t)-f(-t)\over2}$, the even and odd parts
of $f$. For every $0\le t\le s<\alpha$, since $\bmatrix s&-s\\-s&s\endbmatrix
\ge\bmatrix t&-t\\-t&t\endbmatrix\ge0$, we have for $0\le\lambda\le1$
\begin{align*}
&\bmatrix f(\lambda s+(1-\lambda)t)&f(-(\lambda s+(1-\lambda)t))\\
f(-(\lambda s+(1-\lambda)t))&f(\lambda s+(1-\lambda)t)\endbmatrix \\
&\qquad\le\lambda\bmatrix f(s)&f(-s)\\f(-s)&f(s)\endbmatrix
+(1-\lambda)\bmatrix f(t)&f(-t)\\f(-t)&f(t)\endbmatrix\,.
\end{align*}
Multiply $\bmatrix1&\pm1\\0&0\endbmatrix$ from the left and
$\bmatrix1&0\\\pm1&1\endbmatrix$ from the right to get
$$
f(\lambda s+(1-\lambda)t)\pm f(-(\lambda s+(1-\lambda)t)) \\
\le\lambda(f(s)\pm f(-s))+(1-\lambda)(f(t)\pm f(-t)).
$$
Hence $f_0$ and $f_1$ are convex on $[0,\alpha)$. Noting $f=f_0+f_1$, when
$0<s<t<\alpha$, we get
$$
f'_+(s)\le f'_-(t)=(f_0)'_-(t)+(f_1)'_-(t)
\le(f_0)'_+(t)+(f_1)'_+(t)=f'_+(t),
$$
where $f'_-(t)$ is the left-derivative of $f$ at $t$. Thanks to the continuity of
$f'_+$ on $(0,\alpha)$ letting $s\to t$ gives $(f_0)'_-(t)=(f_0)'_+(t)$ and
$(f_1)'_-(t)=(f_1)'_+(t)$, and so $f_0$ and $f_1$ are differentiable on $(0,\alpha)$.
Hence $f=f_0+f_1$ is differentiable on $(-\alpha,\alpha)\setminus\{0\}$ because $f_0$
is even and $f_1$ is odd. For $0<a<\alpha/2$, since the function $f(\cdot+a)$ belongs
to $S_\conv^{(2)}(-\alpha+a,\alpha-a)$, what we have just proven implies that
$f(\cdot+a)$ is differentiable at $-a$ so that $f$ is differentiable at $0$ as well.

Next let us show that $f'\in S_\mono^{(2)}(-\alpha,\alpha)$. Let $A\ge B\ge0$
($A\ne B$) in $M_2(\bR)$ with entries in $(-\alpha,\alpha)$. Write
$A=\bmatrix a_1&a_3\\a_3&a_2\endbmatrix$ and $B=\bmatrix b_1&b_3\\b_3&b_2\endbmatrix$;
then $a_1\ge b_1$, $a_2\ge b_2$ and $(a_3-b_3)^2\le(a_1-b_1)(a_2-b_2)$. Choose
$0\le\delta\le a_1-b_1$ such that $(a_3-b_3)^2=(a_1-b_1-\delta)(a_2-b_2)$, and set
$C:=\bmatrix b_1+\delta&b_3\\b_3&b_2\endbmatrix$. The entries of $C$ are in
$(-\alpha,\alpha)$, and moreover $A\ge C\ge B$ and both $A-C$ and $C-B$ are of at most
rank one. So we may assume that $A-B$ is of rank one so that
$A-B=\bmatrix a&c\\c&b\endbmatrix$ with $a,b\ge0$ and $c^2=ab$. If either $a=0$ or
$b=0$ (hence $c=0$), then $f'[A]\ge f'[B]$ is immediately seen. Now assume $a,b>0$.
Since
$$
f[\lambda A+(1-\lambda)B]\le\lambda f[A]+(1-\lambda)f[B],
\qquad0\le\lambda\le1,
$$
we have for $0<\lambda<1$
\begin{align*}
{f[B+\lambda(A-B)]-f[B]\over\lambda}&\le f[A]-f[B], \\
{f[A+(1-\lambda)(B-A)]-f[A]\over1-\lambda}&\le f[B]-f[A].
\end{align*}
Letting $\lambda\to0$ and $\lambda\to1$ in the above gives
\begin{align*}
(A-B)\circ f'[B]&\le f[A]-f[B], \\
(B-A)\circ f'[A]&\le f[B]-f[A].
\end{align*}
Summing these gives $(A-B)\circ(f'[A]-f'[B])\ge0$. Since
$\bmatrix a^{-1}&c^{-1}\\c^{-1}&b^{-1}\endbmatrix$ (the Schur inverse of $A-B$) is
positive semidefinite, the Schur theorem implies $f'[A]\ge f'[B]$.

To prove the converse, assume that $f$ is differentiable on $(-\alpha,\alpha)$ and
$f'\in\Psi^{(2)}(-\alpha,\alpha)$. Let $A\ge B\ge0$ in $M_2(\bR)$ with entries in
$(-\alpha,\alpha)$. For such $A,B$ we have
\begin{align*}
f\biggl[{A+B\over2}\biggr]-f[B]&={1\over2}\int_0^1(A-B)\circ
f'\biggl[\lambda{A+B\over2}+(1-\lambda)B\biggr]\,d\lambda, \\
{f[A]+f[B]\over2}-f[B]&={f[A]-f[B]\over2}
={1\over2}\int_0^1(A-B)\circ f'[\lambda A+(1-\lambda)B]\,d\lambda.
\end{align*}
Since $f'\in S_\mono^{(2)}(-\alpha,\alpha)$ and
$\lambda A+(1-\lambda)B\ge\lambda{A+B\over2}+(1-\lambda)B\ge0$, we get
$$
f'\biggl[\lambda{A+B\over2}+(1-\lambda)B\biggr]
\le f'[\lambda A+(1-\lambda)B],\qquad0\le\lambda\le1,
$$
and hence $f[{A+B\over2}]\le{f[A]+f[B]\over2}$. Next, since
$$
{k\over2^N}A+\biggl(1-{k\over2^N}\biggr)B
\ge{k-1\over2^N}A+\biggl(1-{k-1\over2^N}\biggr)B,
\qquad k=1,\dots,2^N,\ N\in\bN,
$$
one can easily show by induction that
$$
f\biggl[{k\over2^N}A+\biggl(1-{k\over2^N}\biggr)B\biggr]
\le{k\over2^N}f[A]+\biggl(1-{k\over2^N}\biggr)f[B]
$$
for all $k=0,1,\dots,2^N$ and $N\in\bN$. From the continuity of $A\mapsto f[A]$ we
have $f[\lambda A+(1-\lambda)B]\le\lambda f[A]+(1-\lambda)f[B]$ for all
$0\le\lambda\le1$.
\end{proof}

In the above we characterized functions in the three classes
$S_\pos^{(n)}(-\alpha,\alpha)$, $S_\mono^{(n)}(-\alpha,\alpha)$ and
$S_\conv^{(n)}(-\alpha,\alpha)$ in the case of the first non-trivial order $n=2$.
The following two propositions give necessary conditions for functions
in $S_\pos^{(3)}(-\alpha,\alpha)$ and for those in $S_\mono^{(3)}(-\alpha,\alpha)$,
though complete descriptions of those functions are not known.

\begin{prop}\label{P-2.5}
If $f\in S_\pos^{(3)}(-\alpha,\alpha)$, then it is continuous on $(-\alpha,\alpha)$.
\end{prop}

\begin{proof}
Assume $f\in S_\pos^{(3)}(-\alpha,\alpha)$. Then obviously
$f\in S_\pos^{(2)}(-\alpha,\alpha)$ and by Proposition \ref{P-2.2} it remains to show
the continuity of $f$ on $(-\alpha,0]$. When $f|_{(0,\alpha)}$ is identically zero,
the assertion is obvious from Proposition \ref{P-2.2}. Hence by Lemma \ref{L-2.1} we
may assume that $f>0$ on $(0,\alpha)$. First let us show the right-continuity of $f$
at $0$. For $0<a<\alpha$, since
$$
\bmatrix a&a/\sqrt2&0\\a/\sqrt2&a&a/\sqrt2\\0&a/\sqrt2&a\endbmatrix\ge0,
$$
we get
$$
\bmatrix f(a)&f(a/\sqrt2)&f(0)\\f(a/\sqrt2)&f(a)&f(a/\sqrt2)\\
f(0)&f(a/\sqrt2)&f(a)\endbmatrix\ge0.
$$
Letting $a\searrow0$ gives
$$
\bmatrix f(0+)&f(0+)&f(0)\\f(0+)&f(0+)&f(0+)\\
f(0)&f(0+)&f(0+)\endbmatrix\ge0,
$$
so the determinant is $-f(0+)(f(0+)-f(0))^2\ge0$. Since $0\le f(0)\le f(0+)$ by
Proposition \ref{P-2.2}, we have $f(0+)=f(0)$.

Next let $0\le b<a<\alpha$. Since
$$
\det\bmatrix a&b&t\\b&a&-b\\t&-b&a\endbmatrix
=-(t+a)(at-a^2+2b^2),
$$
we have
$$
\bmatrix a&b&t\\b&a&-b\\t&-b&a\endbmatrix\ge0
\quad\mbox{if }-a\le t\le{a^2-2b^2\over a}.
$$
Here note that $-a<-b<{a^2-2b^2\over a}$ and ${a^2-2b^2\over a}\searrow-b$ as
$a\searrow b$. If $-a\le t\le{a^2-2b^2\over a}$, then
$$
\bmatrix f(a)&f(b)&f(t)\\f(b)&f(a)&f(-b)\\f(t)&f(-b)&f(a)\endbmatrix\ge0
$$
and by taking determinant we get
$$
-f(a)f(t)^2+2f(b)f(-b)f(t)+f(a)^3-f(a)f(b)^2-f(a)f(-b)^2\ge0.
$$
From $f(a)>0$ this gives
\begin{align*}
&{f(b)\over f(a)}f(-b)-
{\sqrt{(f(a)^2-f(b)^2)(f(a)^2-f(-b)^2)}\over f(a)} \\
&\qquad\le f(t)\le {f(b)\over f(a)}f(-b)+
{\sqrt{(f(a)^2-f(b)^2)(f(a)^2-f(-b)^2)}\over f(a)}.
\end{align*}
Therefore,
\begin{align*}
|f(t)-f(-b)|
&\le\bigg|{f(b)\over f(a)}-1\bigg|\cdot|f(-b)|
+\sqrt{f(a)^2-f(b)^2} \\
&\le f(a)-f(b)+\sqrt{f(a)^2-f(b)^2}
\end{align*}
because $|f(-b)|\le f(b)\le f(a)$ by Proposition \ref{P-2.2} and Lemma \ref{L-2.1}.
Since $f(a)\to f(b)$ as $a\searrow b$ (for $b=0$ this was shown above), the above
estimate implies that $f$ is continuous at $-b$ for each $b\in[0,\alpha)$. Hence $f$
is continuous on $(-\alpha,0]$.
\end{proof}
 
\begin{prop}\label{P-2.6}
If $f\in S_\mono^{(3)}(-\alpha,\alpha)$, then it is continuously
differentiable on $(-\alpha,\alpha)$.
\end{prop}

\begin{proof}
Assume $f\in S_\mono^{(3)}(-\alpha,\alpha)$. Since
$f\in S_\mono^{(2)}(-\alpha,\alpha)$, Proposition \ref{P-2.3} implies that $f$ is
differentiable on $(0,\alpha)$, differentiable a.e.\ on $(-\alpha,0)$ and there exist
$g\in\Psi^{(1)}(-\alpha,\alpha)$ and a set $N\subset(-\alpha,0]$ of measure zero
such that $f'(t)=g(t)$ for all $t\in(-\alpha,\alpha)\setminus N$. For every 
$A\in M_3(\bR)$ with entries in $(-\alpha,\alpha)\setminus N$ we have
${f[A+\eps J]-f[A]\over\eps}\ge0$ for all small $\eps>0$. Letting $\eps\searrow0$
gives $g[A]\ge0$. Let us prove that $g$ restricted on $[c,d]\setminus N$ is uniformly
continuous for any closed interval $[c,d]\subset(-\alpha,0)$. One can perform the
argument in the second paragraph of the proof of Proposition \ref{P-2.5} for $g$ in
place of $f$ whenever $0<b<a<\alpha$ and $-b,t\not\in N$. So we see that if
$0<b<a<\alpha$, $-b\not\in N$ and $t\in[-a,{a^2-2b^2\over2}]\setminus N$, then
$$
|g(t)-g(-b)|\le g(a)-g(b)+\sqrt{g(a)^2-g(b)^2}.
$$
Suppose that the asserted uniform continuity is not satisfied. Then for some
$\eps>0$ one can choose $t_k,t'_k\in[c,d]\setminus N$ so that $|t_k-t'_k|\to0$ and
$|g(t_k)-g(t'_k)|\ge\eps$. We may assume $t_k\to t_0$ (also $t'_k\to t_0$) for some
$t_0\in[c,d]$. Since $g$ is continuous on $(0,\alpha)$ and $N$ has measure zero, one
can choose $0<b<a<\alpha$ such that $-b\not\in N$, $-a<t_0<{a^2-2b^2\over2}$ and
$$
g(a)-g(b)+\sqrt{g(a)^2-g(b)^2}<{\eps\over2}.
$$
Since $t_k,t'_k\in[-a,{a^2-2b^2\over2}]\setminus N$ for $k$ large, we get
$$
|g(t_k)-g(t'_k)|\le|g(t_k)-g(-b)|+|g(t'_k)-g(-b)|<\eps,
$$
a contradiction. Hence the uniform continuity of $g$ on $[c,d]\setminus N$ is proven
for any interval $[c,d]\subset(-\alpha,0)$. This implies that
$g|_{(-\alpha,0)\setminus N}$ can extend to a continuous function $\tilde g$ on
$(-\alpha,0)$. Define a function $\tilde g$ on the whole $(-\alpha,\alpha)$ by
$$
\tilde g(t)=\begin{cases}
\tilde g(t) & \text{for $-\alpha<t<0$}, \\
g(0+) & \text{for $t=0$}, \\
g(t) & \text{for $0<t<\alpha$}.
\end{cases}
$$
Note that $g(0+)$ exists since $g$ is nonnegative and non-decreasing on $(0,\alpha)$.
Now we prove that $\tilde g\in S_\pos^{(3)}(-\alpha,\alpha)$. Let
$A=[a_{ij}]\in M_3(\bR)$ with $a_{ij}\in(-\alpha,\alpha)$. One can choose a sequence
$\eps_n\searrow0$ such that $a_{ij}+\eps_n\not\in N$ and $a_{ij}+\eps_n\ne0$ for all
$1\le i,j\le3$ and $n$. Then $g(a_{ij}+\eps_n)\to\tilde g(a_{ij})$ as $n\to\infty$
by definition of $\tilde g$ and $g[A+\eps_nJ]\ge0$ due to $a_{ij}+\eps_n\not\in N$;
so $\tilde g[A]\ge0$ is shown. Hence $\tilde g\in S_\pos^{(3)}(-\alpha,\alpha)$ so
that $\tilde g$ is continuous on $(-\alpha,\alpha)$ by Proposition \ref{P-2.5}. Since
$f'(t)=\tilde g(t)$ a.e.\ on $(-\alpha,\alpha)$, we have
$$
f(t)-f(0)=\int_0^t\tilde g(s)\,ds\quad\mbox{for }-\alpha<t<\alpha.
$$
This implies that $f$ is differentiable on $(-\alpha,\alpha)$ with $f'=\tilde g$.
\end{proof}

\section{Relations among three classes}
\setcounter{equation}{0}

In this section we present some relations among three classes
$S_{\rm pos}^{(n)}(-\alpha,\alpha)$, $S_{\rm mono}^{(n)}(-\alpha,\alpha)$ and
$S_{\rm conv}^{(n)}(-\alpha,\alpha)$ for general $n$.

\begin{prop}\label{P-3.1}
$S_\pos^{(2n)}(-\alpha,\alpha)\subset S_\mono^{(n)}(-\alpha,\alpha)$ and
$S_\mono^{(2n)}(-\alpha,\alpha)\subset S_\conv^{(n)}(-\alpha,\alpha)$ for every
$n\in\bN$.
\end{prop}

\begin{proof}
Assume $f\in S_\pos^{(2n)}(-\alpha,\alpha)$. If $A\ge B\ge0$ in $M_n(\bR)$ with
entries in $(-\alpha,\alpha)$, then
$$
\bmatrix A&B\\B&B\endbmatrix
=\bmatrix A-B&0\\0&0\endbmatrix+\bmatrix B&B\\B&B\endbmatrix\ge0
$$
so that $\bmatrix f[A]&f[B]\\f[B]&f[B]\endbmatrix\ge0$. This implies $f[A]\ge f[B]$
because
$$
\bmatrix f[A]-f[B]&0\\0&0\endbmatrix
=\bmatrix I&-I\\0&0\endbmatrix\bmatrix f[A]&f[B]\\f[B]&f[B]\endbmatrix
\bmatrix I&0\\-I&0\endbmatrix.
$$
Hence $f\in S_\mono^{(n)}(-\alpha,\alpha)$.

Next assume $f\in S_\mono^{(2n)}(-\alpha,\alpha)$ and let $A\ge B\ge0$ in $M_n(\bR)$
with entries in $(-\alpha,\alpha)$. As in the beginning of the proof of Proposition
\ref{P-2.3} by replacing $s,t$ by $A,B$, one can prove that
$f\bigl[{A+B\over2}\bigr]\le{f[A]+f[B]\over2}$. Since $f$ is continuous on
$(-\alpha,\alpha)$, this implies that
$f[\lambda A+(1-\lambda)B]\le\lambda f[A]+(1-\lambda)f[B]$ for all $0\le\lambda\le1$
(see the last of the proof of Proposition \ref{P-2.4}).
\end{proof}

The next theorem extends Proposition \ref{P-2.4} (for $n=2$) and Proposition
\ref{P-2.6} (for $n=3$).

\begin{thm}\label{T-3.2} \
\begin{itemize}
\item[\rm(1)] For every $n\ge2$, $f\in S_\conv^{(n)}(-\alpha,\alpha)$ if and only if
$f$ is differentiable on $(-\alpha,\alpha)$ and $f'\in S_\mono^{(n)}(-\alpha,\alpha)$.
\item[\rm(2)] For every $n\ge3$, $f\in S_\mono^{(n)}(-\alpha,\alpha)$ if and only if
$f$ is differentiable on $(-\alpha,\alpha)$ and $f'\in S_\pos^{(n)}(-\alpha,\alpha)$.
\end{itemize}
\end{thm}

\begin{proof}
(1)\enspace Assume $f\in S_\conv^{(n)}(-\alpha,\alpha)$ with $n\ge2$. Then $f$ is
continuously differentiable on $(-\alpha,\alpha)$ by Proposition \ref{P-2.4}. Let
$A\ge B\ge0$ in $M_n(\bR)$ with entries in $(-\alpha,\alpha)$. Then there are
$A_k\in M_n(\bR)$, $0\le k\le n$, such that
$A=A_0\ge A_1\ge\dots\ge A_{n-1}\ge A_n=B$, all entries of $A_k$'s are in
$(-\alpha,\alpha)$ and $A_{k-1}-A_k$ is of at most rank one for $1\le k\le n$.
In fact, diagonalize $A-B$ as
$A-B=T\,\mathrm{Diag}(\lambda_1,\dots,\lambda_n)\,T^{-1}$ with an orthogonal matrix
$T$ and set
$$
A_k:=B+T\,\mathrm{Diag}(0,\dots,0,\lambda_{k+1},\dots,\lambda_n)\,T^{-1},
\qquad0\le k\le n.
$$
Here we note that all entries of $A_k$'s are in $(-\alpha,\alpha)$ since
$A\ge A_k\ge B\ge0$. Hence we may prove that $f'[A]\ge f'[B]$ if $A\ge B\ge0$ in
$M_n(\bR)$ with entries in $(-\alpha,\alpha)$ and $A-B$ is of rank one. By continuity
of $f'$ we may further assume that $A-B=[a_ia_j]_{1\le i,j\le n}$ with nonzero
$a_1,\dots,a_n\in\bR$ so that $[a_i^{-1}a_j^{-1}]_{1\le i,j\le n}$ is positive
semidefinite. In this situation, the proof of $f'[A]\ge f'[B]$ is same as the second
paragraph of the proof of Proposition \ref{P-2.4}. Moreover, the proof of the
converse is same as the third paragraph of that of Proposition \ref{P-2.4}.

(2)\enspace Assume $f\in S_\mono^{(n)}(-\alpha,\alpha)$ with $n\ge3$. Then $f$ is
differentiable on $(-\alpha,\alpha)$ by Proposition \ref{P-2.6}, and
$f'\in S_\pos^{(n)}(-\alpha,\alpha)$ is seen as in the first part of the proof of
Proposition \ref{P-2.6}. The converse follows as in the second paragraph of the proof
of Proposition \ref{P-2.3}. Here note that $f'$ is continuous on $(-\alpha,\alpha)$
by Proposition \ref{P-2.5}.
\end{proof}

Theorem \ref{T-3.2} further says that for every $n\ge3$,
$f\in S_\conv^{(n)}(-\alpha,\alpha)$ if and only if $f$ is twice differentiable on
$(-\alpha,\alpha)$ and $f''\in S_\pos^{(n)}(-\alpha,\alpha)$.

The next proposition is similar to the obvious fact that if
$f\in S_{\rm mono}^{(n)}(-\alpha,\alpha)$ and $f(0)\ge0$ then
$f\in S_\pos^{(n)}(-\alpha,\alpha)$.

\begin{prop}\label{P-3.3}
For every $n\ge2$, if $f\in S_{\rm conv}^{(n)}(-\alpha,\alpha)$ and $f'(0)\ge0$, then
$f\in S_\mono^{(n)}(-\alpha,\alpha)$.
\end{prop}

\begin{proof}
If $A\ge B\ge0$ in $M_n(\bR)$ with entries in $(-\alpha,\alpha)$, then we have
$$
f[A]-f[B]-(A-B)\circ f'[B]
=(A-B)\circ\int_0^1(f'[\lambda A+(1-\lambda)B]-f'[B])\,d\lambda.
$$
For every $0\le\lambda\le1$, since $\lambda A+(1-\lambda)B\ge B\ge0$, Theorem
\ref{T-3.2}\,(1) implies that $f'[\lambda A+(1-\lambda)B]\ge f'[B]$. Also $f'[B]\ge0$
follows from $f'(0)\ge0$. Hence $f[A]-f[B]\ge(A-B)\circ f'[B]\ge0$ by the Schur
theorem.
\end{proof}

\begin{remark}\label{R-3.4}{\rm
According to \cite[Theorem 1.2]{Ho} and Theorem \ref{T-3.2}, when $n\ge3$ we notice
the following:
\begin{itemize}
\item[(i)] If $f\in S_\pos^{(n)}(-\alpha,\alpha)$, then
$f|_{(0,\alpha)}\in C^{n-3}(0,\alpha)$, $f^{(k)}(x)\ge0$ for all $x\in(0,\alpha)$
and $0\le k\le n-3$, and $f^{(n-3)}$ is non-decreasing and convex on $(0,\alpha)$.
\item[(ii)] If $f\in S_\mono^{(n)}(-\alpha,\alpha)$, then
$f|_{(0,\alpha)}\in C^{n-2}(0,\alpha)$, $f^{(k)}(x)\ge0$ for all $x\in(0,\alpha)$
and $1\le k\le n-2$, and $f^{(n-2)}$ is non-decreasing and convex on $(0,\alpha)$.
\item[(iii)] If $f\in S_\conv^{(n)}(-\alpha,\alpha)$, then
$f|_{(0,\alpha)}\in C^{n-1}(0,\alpha)$, $f^{(k)}(x)\ge0$ for all $x\in(0,\alpha)$
and $2\le k\le n-1$, and $f^{(n-1)}$ is non-decreasing and convex on $(0,\alpha)$.
\end{itemize}
More strongly, it may be expected that if $f\in S_\pos^{(n)}(-\alpha,\alpha)$, then
$f\in C^{n-3}(-\alpha,\alpha)$ and $f|_{(0,\alpha)}\in C^{n-2}(0,\alpha)$. In
particular, it may be conjectured that if $f\in S_\pos^{(3)}(-\alpha,\alpha)$ then
$f|_{(0,\alpha)}\in C^{1}(0,\alpha)$. As will be shown in Section 5 (see
Theorem \ref{T-5.1}), $f(x)=|x|$ is a non-differentiable example in
$S_\pos^{(3)}(-\infty,\infty)$, and $f(x)=(\sign x)x^2$ is in
$S_\mono^{(3)}(-\infty,\infty)$ but it is not twice differentiable. These examples
suggest that the necessary conditions in Propositions \ref{P-2.5} and \ref{P-2.6} are
rather optimal.
}\end{remark}

\section{Characterizations}
\setcounter{equation}{0}

The next theorem characterizes the three classes of S-positive, S-monotone and
S-convex functions on $(-\alpha,\alpha)$. It also shows the explicit differences among
the three notions of S-positivity, S-monotonicity and S-convexity.

\begin{thm}\label{T-4.1}
Let $f$ be a real function on $(-\alpha,\alpha)$,
$0<\alpha\le\infty$. The following statements hold:
\begin{itemize}
\item[\rm(i)] $f$ is S-positive if and only if it is analytic and $f^{(k)}(0)\ge0$
for all $k\ge0$.
\item[\rm(ii)] $f$ is S-monotone if and only if it is analytic and $f^{(k)}(0)\ge0$
for all $k\ge1$.
\item[\rm(iii)] $f$ is S-convex if and only if it is analytic and $f^{(k)}(0)\ge0$
for all $k\ge2$.
\end{itemize}
\end{thm}

For the proof we need the following two lemmas.

\begin{lemma}\label{L-4.2}
Let $f$ be a real function on $(-\alpha,\alpha)$, and $f_0$ and $f_1$ be the even and
odd parts of $f$, i.e., $f_0(x):={f(x)+f(-x)\over2}$ and $f_1(x):={f(x)-f(-x)\over2}$.
Then $f$ is S-convex if and only if so are both $f_0$ and $f_1$.
\end{lemma}

\begin{proof}
Since $f=f_0+f_1$, it is obvious that $f$ is S-convex if so are $f_0$ and $f_1$. To
prove the converse, let $A\ge B\ge0$ in $M_n(\bR)$ with entries in $(-\alpha,\alpha)$
and $0\le\lambda\le1$. If $f$ is S-convex, then we get
\begin{align*}
&\bmatrix f[\lambda A+(1-\lambda)B]
&f[-(\lambda A+(1-\lambda)B)]\\f[-(\lambda A+(1-\lambda)B)]
&f[\lambda A+(1-\lambda)B]\endbmatrix \\
&\qquad\le\lambda\bmatrix f[A]&f[-A]\\f[-A]&f[A]\endbmatrix
+(1-\lambda)\bmatrix f[B]&f[-B]\\f[-B]&f[B]\endbmatrix\,.
\end{align*}
Multiplying $\bmatrix I&\pm I\\0&0\endbmatrix$ from the left and
$\bmatrix I&0\\\pm I&0\endbmatrix$ from the right gives
$$
f[\lambda A+(1-\lambda)B]\pm f[-(\lambda A+(1-\lambda)B)]
\le\lambda(f[A]\pm f[-A])+(1-\lambda)(f[B]\pm f[-B]).
$$
Hence $f_0$ and $f_1$ are S-convex.
\end{proof}

We note that the assertions similar to the above lemma for S-positive and S-monotone
functions are also easy to show.

\begin{lemma}\label{L-4.3}
Let $f$ be an even or odd real function on $(-\alpha,\alpha)$. If $f$ is infinitely
times differentiable on $(-\alpha,\alpha)$ and $f^{(k)}(x)\ge0$ for all
$x\in[0,\alpha)$ and all $k\ge N$ with some $N\in\bN$, then the Taylor expansion
$\sum_{k=0}^\infty(f^{(k)}(0)/k!)x^k$ converges to $f(x)$ for every
$x\in(-\alpha,\alpha)$.
\end{lemma}

\begin{proof}
We may and do assume that $f$ is even and $f^{(k)}(x)\ge0$ for all $x\in[0,\alpha)$
and all $k\ge0$. In fact, we may consider the even function $f'$ when $f$ is odd, and
we obtain the conclusion when the same assertion is proven for $f^{(N)}$ instead of
$f$. For each $0<x<\alpha$ and $m\in\bN$ the Taylor theorem says that
$$
f(x)=\sum_{k=0}^m{f^{(k)}(0)\over k!}x^k+{f^{(m+1)}(\theta x)\over(m+1)!}x^{m+1}
$$
with some $0<\theta<1$. Hence we get $\sum_{k=0}^m(f^{(k)}(0)/k!)x^k\le f(x)$ for all
$0\le x<\alpha$, and so $\sum_{k=0}^\infty(f^{(k)}(0)/k!)x^k$
($=\sum_{k=0}^\infty(f^{(2k)}(0)/(2k)!)x^{2k}$ since $f$ is even)
converges for all $-\alpha<x<\alpha$. Put
$$
g(x):=f(x)-\sum_{k=0}^\infty{f^{(k)}(0)\over k!}x^k
\quad\mbox{for }-\alpha<x<\alpha.
$$
Then $g$ is an infinitely times differentiable even function such that $g(x)\ge0$ for
$0\le x<\alpha$ and $g^{(k)}(0)=0$ for all $k\ge0$. Furthermore, the above argument
applied to $f^{(k)}$ ($k\in\bN$) instead of $f$ shows that
$\sum_{j=0}^\infty(f^{(k+j)}(0)/j!)x^j\le f^{(k)}(x)$ for all $0\le x<\alpha$.
Hence $g^{(k)}(x)\ge0$ for all $0\le x<\alpha$ and all $k\ge0$. To prove that $g$ is
identically zero, suppose that $g(x)>0$ for some $x\in[0,\alpha)$, and let
$\beta:=\inf\{x\in[0,\alpha):g(x)>0\}$. Then $0\le\beta<\alpha$ and $g(x)=0$ for
$-\beta\le x\le\beta$. Define an even function $h$ on $(-\alpha+\beta,\alpha-\beta)$ by
$$
h(x)=\begin{cases}
g(x-\beta) & \text{if $-\alpha+\beta\le x\le0$}, \\
g(x+\beta) & \text{if $0\le x\le\alpha-\beta$}.
\end{cases}
$$
Note that $h$ satisfies the same conditions as $g$, that is, $h^{(k)}(x)\ge0$ for all
$x\in[0,\alpha-\beta)$ and $h^{(k)}(0)=0$ for all $k\ge0$. Let $0<x<(\alpha-\beta)/2$.
For each $m\in\bN$, by the Taylor theorem we have
$$
h(x)={h^{(m)}(\theta x)\over m!}x^m\le{h^{(m)}(x)\over m!}x^m
$$
and
$$
h(2x)=\sum_{k=0}^m{h^{(k)}(x)\over k!}x^k
+{h^{(m+1)}((1+\theta')x)\over(m+1)!}x^{m+1}
\ge\sum_{k=0}^m{h^{(k)}(x)\over k!}x^k
$$
with some $0<\theta,\theta'<1$. The above two inequalities together imply that
$h(x)=0$ for all $0<x<(\alpha-\beta)/2$, contradicting the definition of $h$.
This completes the proof.
\end{proof}

\noindent
{\it Proof of Theorem \ref{T-4.1}.}\enspace
In Example \ref{E-1.3} we saw the ``if" parts of the statements (i)--(iii). Assume
that $f$ is S-convex. By Lemma \ref{L-4.2} we may further assume that $f$ is an even
or odd function. Iterated use of Theorem \ref{T-3.2}\,(1) and Proposition \ref{P-3.1}
implies that $f$ is infinitely times differentiable on $(-\alpha,\alpha)$ and
$f^{(k)}$ is S-convex for all $k\ge0$. In particular, $f^{(k)}$ is convex on
$[0,\alpha)$ for all $k\ge0$, so $f^{(k)}(x)\ge0$ for all $x\in[0,\alpha)$ and all
$k\ge2$. Hence Lemma \ref{L-4.3} proves the ``only if" part of (iii). Assume that $f$
is S-monotone. Then $f$ is S-convex by Proposition \ref{P-3.1}, and $f'(0)\ge0$ is
obvious since $f$ is non-decreasing on $[0,\alpha)$. Hence (ii) is proven. Finally
assume that $f$ is S-positive. Then $f$ is S-monotone by Proposition \ref{P-3.1} and
$f(0)\ge0$ is obvious. Hence (i) holds.
\qed

\bigskip
Since the arguments in Example \ref{E-1.3} are valid for functions $z^k$ on $\bC$ and
complex matrices $A\ge B\ge0$ as well, Theorem \ref{T-4.1} yields the following:

\begin{cor}\label{C-4.4}
If $f:(-\alpha,\alpha)\to\bR$ is S-positive, then $f$ has a complex analytic
continuation $\tilde f$ on $\{z\in\bC:|z|<\alpha\}$ and $\tilde f$ is S-positive in
the sense that $[\tilde f(a_{ij})]\ge0$ for all $A=[a_{ij}]\ge0$ in $M_n(\bC)$ of any
$n$ with $|a_{ij}|<\alpha$ for all $i,j$. The similar statements are valid also for
an S-monotone function or an S-convex function.
\end{cor}

The following are typical examples:
\begin{itemize}
\item $f(x):=e^x$ is S-positive on $(-\infty,\infty)$,
\item $f(x):=-\log(1-x)=\sum_{k=1}^\infty(1/k)x^k$ convergent for $|x|<1$ is
S-positive on $(-1,1)$,
\item for $0<p<1$, $f(x):=-(1-x)^p=-1+\sum_{k=1}^\infty(-1)^{k-1}{p\choose k}x^k$
convergent for $|x|\le1$ is S-monotone on $(-1,1)$.
\end{itemize}

\section{Examples of fractional power functions}
\setcounter{equation}{0}

For $p>0$ define an even function $\phi_p$ and an odd function $\psi_p$ on $\bR$ by
$$
\phi_p(x):=|x|^p,\quad\psi_p(x):=(\sign x)|x|^p\quad\mbox{for }x\in\bR.
$$
Also set $\phi_0(x):=1$ and $\psi_0(x):=\sign x$, i.e., $\psi_0(x):=-1,0,1$ if
$x<0$, $x=0$, $x>0$, respectively. The next theorem extends
\cite[Theorems 2.2 and 2.4]{FH}.

\begin{thm}\label{T-5.1} \
\begin{itemize}
\item[\rm(i)] If $n\ge2$ and $p\ge n-2$, then
$\phi_p,\psi_p\in S_\pos^{(n)}(-\infty,\infty)$.
\item[\rm(ii)] If $n\ge1$ and $p\ge n-1$, then
$\phi_p,\psi_p\in S_\mono^{(n)}(-\infty,\infty)$.
\item[\rm(iii)] If $n\ge1$ and $p\ge n$, then
$\phi_p,\psi_p\in S_\conv^{(n)}(-\infty,\infty)$.
\end{itemize}
\end{thm}

\begin{proof}
(i)\enspace Prove by induction on $n$. When $n=2$, the assertion is immediately seen.
Also, when $n=3$ and $p=1$, the result for $\phi_1(x)=|x|$ is well known, and that
for $\psi_1(x)=x$ is trivial. Next assume that the assertion holds for $n\ge2$, and
assume $p\ge n-1$ with $p>1$. (When $p=1$ and so $n=2$, the assertion for $n+1=3$
holds as mentioned above.) For $p>1$ note that $\phi_p$ and $\psi_p$ are
differentiable as
$$
\phi_p'(x)=p\psi_{p-1}(x),\quad\psi_p'(x)=p\phi_{p-1}(x)
\quad\mbox{for }x\in\bR.
$$
Now let us proceed as in the proof of \cite[Theorem 2.2]{FH}. Let $A=[a_{ij}]\ge0$
in $M_{n+1}(\bR)$. Let
$\xi:=(a_{1,n+1},a_{2,n+1},\dots,a_{n+1,n+1})^t/\sqrt{a_{n+1,n+1}}$ if
$a_{n+1,n+1}>0$, and $\xi$ be the zero vector if $a_{n+1,n+1}=0$. Then
$A-\xi\xi^t\ge0$ and we have
\begin{align*}
\phi_p[A]&=\phi_p[\xi\xi^t]+p\int_0^1(A-\xi\xi^t)
\circ\psi_{p-1}[\lambda A+(1-\lambda)\xi\xi^t]\,d\lambda, \\
\psi_p[A]&=\psi_p[\xi\xi^t]+p\int_0^1(A-\xi\xi^t)
\circ\phi_{p-1}[\lambda A+(1-\lambda)\xi\xi^t]\,d\lambda.
\end{align*}
Since the last row and column of $A-\xi\xi^t$ are zero and since $p-1\ge n-2$, the
induction hypothesis and the Schur theorem imply that
\begin{align*}
&\int_0^1(A-\xi\xi^t)\circ
\psi_{p-1}[\lambda A+(1-\lambda)\xi\xi^t]\,d\lambda\ge0, \\
&\int_0^1(A-\xi\xi^t)\circ
\phi_{p-1}[\lambda A+(1-\lambda)\xi\xi^t]\,d\lambda\ge0.
\end{align*}
Furthermore, $\phi_p[\xi\xi^t]\ge0$ and $\psi_p[\xi\xi^t]\ge0$ are immediately seen.
Hence $\phi_p[A]\ge0$ and $\psi_p[A]\ge0$ so that the assertion for $n+1$ is proven.

\medskip
(ii)\enspace The assertion is trivial for $n=1$. When $n=2$ and $p=1$, the result for
$\phi_1(x)=|x|$ is easy to check, and that for $\psi_1(x)=x$ is trivial. Now let
$n\ge2$, and assume $p\ge n-1$ with $p>1$. Let $A\ge B\ge0$ in $M_n(\bR)$. We have
\begin{align*}
\phi_p[A]-\phi_p[B]&=p\int_0^1(A-B)\circ
\psi_{p-1}[\lambda A+(1-\lambda)B]\,d\lambda, \\
\psi_p[A]-\psi_p[B]&=p\int_0^1(A-B)\circ
\phi_{p-1}[\lambda A+(1-\lambda)B]\,d\lambda.
\end{align*}
Since $p-1\ge n-2$, the above (i) gives $\psi_{p-1}[\lambda A+(1-\lambda)B]\ge0$ and
$\phi_{p-1}[\lambda A+(1-\lambda)B]\ge0$. Hence we obtain $\phi_p[A]\ge\phi_p[B]$ and
$\psi_p[A]\ge\psi_p[B]$.

\medskip
(iii)\enspace It suffices to prove the inequalities for $\lambda=1/2$ (see the last of
the proof of Proposition \ref{P-2.4}). When $n=1$ and $p=1$, the result is trivial.
Assume $p\ge n$ with $p>1$, and let $A\ge B\ge0$ in $M_n(\bR)$. We have
\begin{align*}
\phi_p\biggl[{A+B\over2}\biggr]-\phi_p[B]&={p\over2}\int_0^1(A-B)\circ
\psi_{p-1}\biggl[\lambda{A+B\over2}+(1-\lambda)B\biggr]\,d\lambda, \\
{\phi_p[A]+\phi_p[B]\over2}-\phi_p[B]&={\phi_p[A]-\phi_p[B]\over2}
={p\over2}\int_0^1(A-B)\circ\psi_{p-1}[\lambda A+(1-\lambda)B]\,d\lambda.
\end{align*}
Since $p-1\ge n-1$ and $\lambda A+(1-\lambda)B\ge\lambda{A+B\over2}+(1-\lambda)B\ge0$,
the above (ii) implies that
$$
\phi_{p-1}\biggl[\lambda{A+B\over2}+(1-\lambda)B\biggr]
\le\phi_{p-1}[\lambda A+(1-\lambda)B],
$$
and hence $\phi_p\bigl[{A+B\over2}\bigr]\le{\phi_p[A]+\phi_p[B]\over2}$. The proof of
$\psi_p\bigl[{A+B\over2}\bigr]\le{\psi_p[A]+\psi_p[B]\over2}$ is similar.
\end{proof}

According to \cite[Theorems 2.2 and 2.4]{FH}, the conditions $p\ge n-2$ in (i) and
$p\ge n-1$ in (ii) of the above theorem are sharp for both $\phi_p$ and $\psi_p$. The
next lemma says that this is the case also for the condition $p\ge n$ in (iii).

\begin{lemma}\label{L-5.2}
If $n\in\bN$ and $0<p<n$ and if $p$ is not an integer, then there exist
$[a_{ij}]\ge[b_{ij}]\ge0$ in $M_n(\bR)$ such that $a_{ij},b_{ij}>0$ for all $i,j$ and
$$
\biggl[\biggl({a_{ij}+b_{ij}\over2}\biggr)^p\biggr]\not\le
\biggl[{a_{ij}^p+b_{ij}^p\over2}\biggr].
$$
\end{lemma}

\begin{proof}
We use the example in the proof of \cite[Theorem 2.2]{FH} and slightly modify the
argument there. Assume that $0<p<n$ is not an integer. Put $A:=[1+ij]_{1\le i,j\le n}$
and $B:=J$. Then $A\ge B\ge0$. Choose a real $n$-vector $\eta=(\eta_1,\dots,\eta_n)^t$
which is orthogonal to $(1^k,2^k,\dots,n^k)^t$ for $2\le k\le[p]+1$ and
$\sum_{i=1}^ni^{[p]+2}\eta_i=1$. Let $A_t:=tA+(1-t)B=[1+tij]$ and moreover
$$
f(t):=\<\phi_p[A_t]\eta,\eta\>
=\sum_{i,j=1}^n(1+tij)^p\eta_i\eta_j
$$
for $-n^{-2}<t<n^{-2}$. The Taylor expansion of $f(t)$ is given as follows:
\begin{align*}
f(t)&=\sum_{i,j=1}^n\sum_{k=0}^\infty
{p\choose k}t^ki^kj^k\eta_i\eta_j
=\sum_{k=0}^\infty{p\choose k}
\Biggl(\sum_{i=1}^ni^k\eta_i\Biggr)^2t^k \\
&=\Biggl(\sum_{i=1}^n\eta_i\Biggr)^2
+p\Biggl(\sum_{i=1}^ni\eta_i\Biggr)^2t
+{p\choose[p]+2}t^{[p]+2}
+\sum_{k=[p]+3}^\infty{p\choose k}
\Biggl(\sum_{i=1}^ni^k\eta_i\Biggr)^2t^k.
\end{align*}
Hence
$$
f''(t)={p\choose[p]+2}t^{[p]}+O(t^{[p]+1})
\quad\mbox{as }\,t\to0.
$$
Since ${p\choose[p]+2}<0$, we get $f''(t)<0$ for $t>0$ sufficiently small. This means
that $f(t)$ is not convex on $[0,\delta]$ for some small $\delta>0$. So there are
$s,t\in[0,\delta]$ such that $f({s+t\over2})>{f(s)+f(t)\over2}$, which implies
$\phi_p\bigl[{A_s+A_t\over2}\bigr]\not\le{\phi_p[A_s]+\phi_p[A_t]\over2}$.
\end{proof}

By Theorem \ref{T-5.1} and Example \ref{E-1.3} the following is immediately seen:
Let $a_0,a_1,\dots$ be nonnegative real numbers and $\mu,\nu$ be positive measures
on $[0,\infty)$ with $\int_0^\infty\alpha^p\,d\mu(p)<+\infty$ and
$\int_0^\infty\alpha^p\,d\nu(p)<+\infty$. For $m=0,1,\dots$ define
$$
f_m(x):=\sum_{k=0}^\infty a_kx^k+\int_m^\infty x^p\,d\mu(p)
+\int_m^\infty(\sign x)|x|^p\,d\nu(p)
\quad\mbox{for }-\alpha<x<\alpha.
$$
Then $f_{n-2}\in S_\pos^{(n)}(-\infty,\infty)$ for any $n\ge2$ and
$f_{n-1}\in S_\mono^{(n)}(-\infty,\infty)$, $f_n\in S_\conv^{(n)}(-\infty,\infty)$ for
any $n\ge1$.

\section{Weak majorizations}
\setcounter{equation}{0}

In this section we give weak majorizations and unitarily invariant norm inequalities
involving entrywise matrix functions. For a Hermitian $n\times n$ matrix $A$ let
$\lambda(A)=(\lambda_1(A),\dots,\lambda_n(A))$ be the eigenvalues of $A$ in
decreasing order, and $s(A)=(s_1(A),\dots,s_n(A))$ be the singular values of $A$ in
decreasing order. For real $n$-vectors $\mathbf{a}=(a_1,\dots,a_n)$ and
$\mathbf{b}=(b_1,\dots,b_n)$ the {\it weak majorization} $\mathbf{a}\prec_w\mathbf{b}$
means that
$$
\sum_{i=1}^ka_{[i]}\le\sum_{i=1}^kb_{[i]}
\quad\mbox{for }1\le k\le n,
$$
where $(a_{[1]},\dots,a_{[n]})$ is the decreasing rearrangement of the coordinates of
$\mathbf{a}$ and similarly for $\mathbf{b}$. The {\it majorization}
$\mathbf{a}\prec\mathbf{b}$ is referred to if in addition equality holds for $k=n$ in
the above (see \cite{Bh,MO} for details on (weak) majorization theory for vectors and
matrices). We write $\mathbf{a}\circ\mathbf{b}$ for the coordinatewise product
$(a_1b_1,\dots,a_nb_n)$ (i.e., the Schur product when regarded as diagonal matrices).

Let $f$ be a real differentiable function on an interval $(\beta,\gamma)$. The
{\it divided difference} of $f$ is the function $f^{[1]}(a,b)$ on $(\beta,\gamma)^2$
defined by
$$
f^{[1]}(a,b):=\begin{cases}
{f(a)-f(b)\over a-b} & \text{if $a\ne b$}, \\
f'(a) & \text{if $a=b$}.
\end{cases}
$$
Moreover, the second divided difference $f^{[2]}(a,b,c)$ on $(\beta,\gamma)^3$ is
defined by
$$
f^{[2]}(a,b,c):={f^{[1]}(a,b)-f^{[1]}(b,c)\over a-c}
$$
under the assumption of $f$ being twice differentiable. In particular,
$$
f^{[2]}(a,b,b)={f(a)-f(b)-f'(b)(a-b)\over(a-b)^2},
\quad f^{[2]}(a,a,a)={1\over2}f''(a).
$$

The next theorem extends \cite[Corollary 1]{Bo}. Here it should be noted that
\cite[Corollary 1]{Bo} is not true without the assumption $f(0)=0$.

\begin{thm}\label{T-6.1}
Assume that $n\ge2$ and $f\in S_\mono^{(n)}(-\alpha,\alpha)$. For every $A\ge0$ in
$M_n(\bR)$ with $\|A\|\le\alpha$,
$$
\lambda(f[A]-f(0)J)\prec_w\lambda(f(A)-f(0)I).
$$
\end{thm}

To prove this, we need an elementary lemma, whose proof is given since we cannot
find a suitable reference.

\begin{lemma}\label{L-6.2}
Let $f$ be a continuous function on $[0,\alpha)$ which is continuously differentiable
on $(0,\alpha)$. If $f(0)\ge0$ and $f'$ is convex on $(0,\alpha)$, then $f(x)/x$ is
convex on $(0,\alpha)$.
\end{lemma}

\begin{proof}
First assume in addition that $f$ has the third derivative on $[0,\alpha)$. For every
$x\in(0,\alpha)$ the Taylor theorem implies that
$$
0\le f(0)=f(x)+f'(x)(0-x)+{f''(x)\over2}(0-x)^2+{f'''(\theta x)\over6}(0-x)^3
$$
for some $0<\theta<1$. Hence
$$
\biggl({f(x)\over x}\biggr)''
={2\over x^3}\biggl(f(x)-f'(x)x+{f''(x)\over2}x^2\biggr)
\ge{f'''(\theta x)\over3}\ge0
$$
so that $f(x)/x$ is convex on $(0,\alpha)$. To prove the lemma without the existence
of the third derivative, let $\phi$ be a smooth function on $\bR$ supported on
$[-1,-1/2]$ such that $\phi(x)\ge0$ and $\int_{-1}^{-1/2}\phi(x)\,dx=1$. For $\eps>0$
set $\phi_\eps(x):=\eps^{-1}\phi(\eps^{-1}x)$, supported on $[-\eps,-\eps/2]$, and
$f_\eps(x):=\int_{-\eps}^{-\eps/2}f(x-t)\phi_\eps(t)\,dt$ for $0\le x<\alpha-\eps$.
Then one can easily see that $f_\eps(x)$ is smooth on $[0,\alpha-\eps)$ and
$f(x)=\lim_{\eps\searrow0}f_\eps(x)$ for all $0\le x<\alpha$. Since
$f_\eps'(x)=\int_{-\eps}^{-\eps/2}f'(x-t)\phi_\eps(t)\,dt$ is convex on
$(0,\alpha-\eps)$, the above case implies that $(f_\eps(x)-f_\eps(0))/x$ is convex on
$(0,\alpha-\eps)$. Since
$$
{f(x)\over x}
=\lim_{\eps\searrow0}{f_\eps(x)-f_\eps(0)\over x}+{f(0)\over x}
\quad\mbox{for }0<x<\alpha,
$$
the conclusion follows.
\end{proof}

\noindent
{\it Proof of Theorem \ref{T-6.1}.}\enspace
Let $A\ge0$ in $M_n(\bR)$ with $\|A\|<\alpha$. First note that $f[A]$ as well as
$f(A)$ can be defined, $f[A]\ge f(0)J$ and $f(A)\ge f(0)I$. We may assume $f(0)=0$
and prove $\lambda(f[A])\prec_w\lambda(f(A))$. Let
$d(A)=\bigl(d_1(A),\dots,d_n(A)\bigr)$ be the diagonal entries of $A$ in decreasing
order.

We begin with the case $n=2$. So assume that
$f\in S_\mono^{(2)}(-\alpha,\alpha)$ (with $f(0)=0$) and
$A=\bmatrix a&c\\c&b\endbmatrix\ge0$ in $M_2(\bR)$ with $\|A\|<\alpha$. It suffices
to show $\|f[A]\|\le\|f(A)\|$ and $\Tr f[A]\le\Tr f(A)$. Let $s:=\|A\|>0$; then
$\|f(A)\|=f(s)$ by Proposition \ref{P-2.3}. For any unit vector
$\xi=\bmatrix\xi_1\\\xi_2\endbmatrix\in\bC^2$ we get
\begin{align*}
|\<f[A]\xi,\xi\>|
&\le f(a)|\xi_1|^2+2|f(c)|\,|\xi_1|\,|\xi_2|+f(b)|\xi_2|^2 \\
&\le f(a)|\xi_1|^2+2f(|c|)\,|\xi_1|\,|\xi_2|+f(b)|\xi_2|^2 
\end{align*}
by Proposition \ref{P-2.3}. Furthermore, since $f(x)\le(f(s)/s)x$ for $0\le x\le s$
thanks to the convexity of $f$ on $[0,\alpha)$ (see the proof of Proposition
\ref{P-2.3}), we have
\begin{align*}
|\<f[A]\xi,\xi\>|
&\le{f(s)\over s}(a|\xi_1|^2+2|c|\,|\xi_1|\,|\xi_2|+b|\xi_2|^2) \\
&\le{f(s)\over s}\bigg\|\bmatrix a&|c|\\|c|&b\endbmatrix\bigg\|
={f(s)\over s}\|A\|=f(s)
\end{align*}
so that $\|f[A]\|\le f(s)=\|f(A)\|$. From the Schur majorization $d(A)\prec\lambda(A)$
and the convexity of $f$ on $[0,\alpha)$ we also get $f(d(A))\prec_wf(\lambda(A))$,
which implies that $\Tr f[A]\le\Tr f(A)$. Hence the case $n=2$ is shown.

Next assume that $n\ge3$ and $f\in S_\mono^{(n)}(-\alpha,\alpha)$ (with $f(0)=0$).
By Theorem \ref{T-3.2}\,(2), $f$ is differentiable on $(-\alpha,\alpha)$ and
$f'\in S_\pos^{(n)}(-\alpha,\alpha)$. It is known \cite[Theorem 1.2]{Ho} (or Remark
\ref{R-3.4}) that $f'$ is nonnegative non-decreasing and convex on $[0,\alpha)$. Set
$g(x):=f(x)/x$ for $x\in(0,\alpha)$; then $g$ is non-decreasing and convex on
$(0,\alpha)$ by Lemma \ref{L-6.2}. By continuity we may assume that $A$ is positive
and invertible. Since $f[A]=A\circ\int_0^1f'[tA]\,dt$ and $\int_0^1f'[tA]\,dt$ as
well as $A$ is positive semidefinite, it follows from the majorization result in
\cite[Theorem 3\,(i)]{BS} that
$$
\lambda(f[A])\prec_w\lambda(A)\circ d\biggl(\int_0^1f'[tA]\,dt\biggr)
=\lambda(A)\circ\int_0^1f'(td(A))\,dt=\lambda(A)\circ g(d(A))
$$
thanks to the non-decreasingness of $f'$. Since $d(A)\prec\lambda(A)$ as already
mentioned, we get $g(d(A))\prec_w g(\lambda(A))=\lambda(g(A))$ thanks to the
convexity and the non-decreasingness of $g$. Therefore,
$$
\lambda(f[A])\prec_w\lambda(A)\circ\lambda(g(A))=\lambda(Ag(A))=\lambda(f(A)),$$
completing the proof.\qed

\bigskip
The weak majorization in Theorem \ref{T-6.1} gives the norm inequality
$$
|||f[A]-f(0)J|||\le|||f(A)-f(0)I|||
$$
for every unitarily invariant norm $|||\cdot|||$. For example, we notice by Theorem
\ref{T-5.1}\,(ii) that if $n\ge2$ and $p\ge n-1$, then
$$
|||\phi_p[A]|||\le|||A^p|||,\quad|||\psi_p[A]|||\le|||A^p|||
$$
for all $A\ge0$ in $M_n(\bR)$ and every unitarily invariant norm.

\begin{thm}\label{T-6.3}
Assume that $n\ge2$ and $f\in S_\conv^{(n)}(-\alpha,\alpha)$ with $f'(0)\ge0$. For
every $A,B\ge0$ in $M_n(\bR)$ with $\|A\|,\|B\|<\alpha$,
$$
s(f[A]-f[B])\prec_wf^{[1]}(\lambda(A),\lambda(B))\circ s(A-B),
$$
where
$$
f^{[1]}(\lambda(A),\lambda(B)):=
\bigl(f^{[1]}(\lambda_1(A),\lambda_1(B)),
\dots,f^{[1]}(\lambda_n(A),\lambda_n(B))\bigr).
$$
\end{thm}

\begin{proof}
Set $g:=f'-f'(0)$, which is in $S_\mono^{(n)}(-\alpha,\alpha)$ by Theorem
\ref{T-3.2}\,(1). Since
$$
f[A]-f[B]=f'(0)(A-B)+(A-B)\circ\int_0^1g[tA+(1-t)B]\,dt,
$$
we have
$$
s(f[A]-f[B])\prec_w
f'(0)s(A-B)+s(A-B)\circ\lambda\biggl(\int_0^1g[tA+(1-t)B]\,dt\biggr)
$$
by the Ky Fan majorization theorem (\cite[p.\,243]{MO}, \cite[(II.18)]{Bh}) and by
the majorization result \cite[Lemma 1]{HJ} (independently \cite[Lemma 1]{O}), noting
that $g[tA+(1-t)B]\ge0$ for all $0\le t\le1$. The Ky Fan majorization theorem again
gives
$$
\lambda\biggl(\int_0^1g[tA+(1-t)B]\,dt\biggr)
\prec\int_0^1\lambda\bigl(g[tA+(1-t)B]\bigr)\,dt.
$$
Furthermore, Theorem \ref{T-6.1} implies that
$$
\lambda\bigl(g[tA+(1-t)B]\bigr)\prec_w\lambda\bigl(g(tA+(1-t)B)\bigr)
=g\bigl(\lambda(tA+(1-t)B)\bigr)
$$
for all $0\le t\le1$. Since $\lambda(tA+(1-t)B)\prec t\lambda(A)+(1-t)\lambda(B)$ and
$g$ is convex on $[0,\alpha)$, we get
$g\bigl(\lambda(tA+(1-t)B)\bigr)\prec_w g\bigl(t\lambda(A)+(1-t)\lambda(B)\bigr)$
so that
$$
\lambda\biggl(\int_0^1g[tA+(1-t)B]\,dt\biggr)
\prec_w\int_0^1g\bigl(t\lambda(A)+(1-t)\lambda(B)\bigr)\,dt.
$$
Here recall the simple fact that if $\mathbf{a}$, $\mathbf{b}$ and $\mathbf{c}$ are
$n$-vectors with nonnegative coordinates in decreasing order, then
$\mathbf{b}\prec_w\mathbf{c}$ implies
$\mathbf{a}\circ\mathbf{b}\prec_w\mathbf{a}\circ\mathbf{c}$.
Hence we obtain
\begin{align*}
s(f[A]-f[B])&\prec_wf'(0)s(A-B)+s(A-B)\circ
\int_0^1g\bigl(t\lambda(A)+(1-t)\lambda(B)\bigr)\,dt \\
&=s(A-B)\circ\int_0^1f'\bigl(t\lambda(A)+(1-t)\lambda(B)\bigr)\,dt \\
&=s(A-B)\circ f^{[1]}(\lambda(A),\lambda(B)),
\end{align*}
as desired.
\end{proof}

Proposition \ref{P-3.3} says that the assumption of Theorem \ref{T-6.1} is
weaker than that of Theorem \ref{T-6.3}. Also, the weak majorization in Theorem
\ref{T-6.1} is the particular case of that of Theorem \ref{T-6.3} when $B=0$. In fact,
notice that $f^{[1]}(\lambda(A),(0,\dots,0))\circ\lambda(A)=\lambda(f(A)-f(0)I)$ when
$f$ and $A$ are as in Theorem \ref{T-6.3}.

\begin{remark}\label{R-6.4}{\rm \
\begin{itemize}
\item[(1)] For $0<p<1$ we have $\phi_p\in S_\pos^{(2)}(-\infty,\infty)$ by Theorem
\ref{T-5.1}\,(i). When $A=\bmatrix1/2&1/2\\1/2&1/2\endbmatrix$ we compute
$\phi_p[A]=2^{1-p}A$ and $A^p=A$ so that $\|\phi_p[A]\|=2^{1-p}>1=\|A^p\|$. Hence
Theorem \ref{T-6.1} is not valid if the assumption $f\in S_\mono^{(n)}(-\alpha,\alpha)$
is weakened to $f\in S_\pos^{(n)}(-\alpha,\alpha)$.
\item[(2)] Let $A=\bmatrix1&1\\1&1\endbmatrix$ and $B=\bmatrix1&-1\\-1&1\endbmatrix$.
For $p>0$ we have $\psi_p[A]=A$, $\psi_p[B]=B$, $s(A-B)=(2,2)$,
$\lambda(A)=\lambda(B)=(2,0)$ and $(x^p)^{[1]}(\lambda(A),\lambda(B))=(p2^{p-1},0)$.
If the weak majorization in Theorem \ref{T-6.3} holds for $\psi_p$, then we must have
$4\le2p2^{p-1}$ and so $2^{2-p}\le p$, which gives $p\ge1.4\cdots$. Hence we notice
by Theorem \ref{T-5.1}\,(ii) that Theorem \ref{T-6.3} is not valid if the assumption
$f\in S_\conv^{(n)}(-\alpha,\alpha)$ is weakened to
$f\in S_\mono^{(n)}(-\alpha,\alpha)$.
\end{itemize}
}\end{remark}

\begin{prop}\label{P-6.5}
Assume that $n\ge2$ and $f$ is differentiable on $(-\alpha,\alpha)$ with
$f'\in S_\conv^{(n)}(-\alpha,\alpha)$ and $f''(0)\ge0$. For every $A,B\ge0$ in
$M_n(\bR)$ with $\|A\|,\|B\|<\alpha$,
\begin{align*}
s(f[A]-f[B]-(A-B)\circ f'[B])
&\prec_wf^{[2]}(\lambda(A),\lambda(B),\lambda(B))\circ s((A-B)\circ(A-B)) \\
&\prec_wf^{[2]}(\lambda(A),\lambda(B),\lambda(B))\circ s((A-B)^2),
\end{align*}
where
$$
f^{[1]}(\lambda(A),\lambda(B),\lambda(B)):=
\bigl(f^{[1]}(\lambda_1(A),\lambda_1(B),\lambda_1(B)),
\dots,f^{[1]}(\lambda_n(A),\lambda_n(B),\lambda_n(B))\bigr).
$$
\end{prop}

\begin{proof}
Set $g(z):=f''(z)-f''(0)$, which is in $S_\mono^{(n)}(-\alpha,\alpha)$. Since
\begin{align*}
&f[A]-f[B]-(A-B)\circ f'[B] \\
&\qquad=f''(0)(A-B)\circ(A-B)
+(A-B)\circ(A-B)\circ\int_0^1dt\int_0^tg[uA+(1-u)B)]\,du,
\end{align*}
we have
\begin{align*}
&s(f[A]-f[B]-(A-B)\circ f'[B]) \\
&\qquad\prec_ws((A-B)\circ(A-B))\circ
\int_0^1dt\int_0^t f''\bigl(u\lambda(A)+(1-u)\lambda(B)\bigr)\,du \\
&\qquad=s((A-B)\circ(A-B))\circ f^{[2]}(\lambda(A),\lambda(B),\lambda(B))
\end{align*}
similarly to the proof of Theorem \ref{T-6.3}. The second weak majorization follows
from $s((A-B)\circ(A-B))\prec_ws((A-B^2)$.
\end{proof}

The following propositions are more weak majorizations of similar vein.

\begin{prop}\label{P-6.6}
Assume that $n\ge3$ and $f\in S_{\rm mono}^{(n)}(-\alpha,\alpha)$. For every $A,B\ge0$
in $M_n(\bR)$ with entries in $(-\alpha,\alpha)$,
$$
s(f[A]-f[B])\prec_w
\biggl(\max_{1\le i\le n}f^{[1]}(a_{ii},b_{ii})\biggr)s(A-B).
$$
\end{prop}

\begin{proof}
By Theorem \ref{T-3.2}\,(2), $f$ is differentiable with
$f'\in S_{\rm pos}^{(n)}(-\alpha,\alpha)$. Applying \cite[Theorem 3]{AHJ} to
$$
f[A]-f[B]=(A-B)\circ\int_0^1f'[B+t(A-B)]\,dt
$$
gives
\begin{align*}
s(f[A]-f[B])&\prec_w\biggl(\max_{1\le i\le n}\int_0^1
f'(b_{ii}+t(a_{ii}-b_{ii}))\,d\lambda\biggr)s(A-B) \\
&=\biggl(\max_{1\le i\le n}f^{[1]}(a_{ii},b_{ii})\biggr)s(A-B),
\end{align*}
as desired.
\end{proof}

It is clear from Proposition \ref{P-2.3} and the above proof that Proposition
\ref{P-6.6} holds for $n=2$ as well whenever $f$ is continuously differentiable. We
remark that the two weak majorizations in Theorem \ref{T-6.3} and Proposition
\ref{P-6.6} are not comparable in general, that is, the right-hand sides of those are
not generally comparable in weak majorization. For example, when $f(x)=x^2$,
$A=\bmatrix1&1&0\\1&1&0\\0&0&0\endbmatrix$ and
$B=\bmatrix2&0&0\\0&0&0\\0&0&0\endbmatrix$, the right-hand side in Theorem \ref{T-6.3}
is $(4\sqrt2,0,0)$ and that of Proposition \ref{P-6.6} is $(3\sqrt2,3\sqrt2,0)$, which
are not comparable.

\begin{prop}\label{P-6.7}
Assume that $n\ge3$ and $f\in S_{\rm conv}^{(n)}(-\alpha,\alpha)$. For every $A,B\ge0$
in $M_n(\bR)$ with entries in $(-\alpha,\alpha)$,
\begin{align*}
s(f[A]-f[B]-(A-B)\circ f'[B])
&\prec_w\biggl(\max_{1\le i\le n}f^{[2]}(a_{ii},b_{ii},b_{ii})\biggr)
s((A-B)\circ(A-B)) \\
&\prec_w\biggl(\max_{1\le i\le n}f^{[2]}(a_{ii},b_{ii},b_{ii})\biggr)
s((A-B)^2).
\end{align*}
\end{prop}

\begin{proof}
By Theorem \ref{T-3.2}, $f$ is twice differentiable with
$f''\in S_{\rm pos}^{(n)}(-\alpha,\alpha)$. Then the proof is similar to those of
Propositions \ref{P-6.5} and \ref{P-6.6}.
\end{proof}

Corresponding to the weak majorizations obtained above, we get unitarily invariant
norm inequalities for entrywise matrix functions. For instance, when $f$, $A$ and $B$
are as in Proposition \ref{P-6.6}, we have
$$
|||f[A]-f[B]|||
\le\Bigl(\max_if^{[1]}(a_{ii},b_{ii})\Bigr)|||A-B|||
\le f^{[1]}\Bigl(\max_ia_{ii},\max_ib_{ii}\Bigr)|||A-B|||
$$
for any unitarily invariant norm $|||\cdot|||$. In particular, we have norm
inequalities as above for the functions $\phi_p,\psi_p$ in Section 5.

Assume that $f$ is S-monotone on $(-\alpha,\alpha)$. By Corollary \ref{C-4.4}, $f$
extends to a complex analytic function, denoted by the same $f$. Then the weak
majorization in Theorem \ref{T-6.3} holds more generally for every $A,B\ge0$ in
$M_n(\bC)$ of any $n$ with $\|A\|,\|B\|<\alpha$, and that in Proposition \ref{P-6.6}
holds for every $A,B\in M_n(\bC)$ of any $n$ with $|a_{ij}|,|b_{ij}|<\alpha$. The
proofs of those generalizations are same as given above. When $f$ is S-convex on
$(-\alpha,\alpha)$, similar generalizations work for the weak majorizations in
Propositions \ref{P-6.5} and \ref{P-6.7}. We thus have the corresponding norm
inequalities, for example, for the functions mentioned at the end of Section 4.

\section*{Acknowledgments}
In 2001 Xingzhi Zhan proposed the author the conjecture that FitzGerald and Horn's
results \cite{FH} can extend to the function $|z|^pe^{i\theta}$, $p>0$, of complex
variable $z=|z|e^{i\theta}$. The first motivation of this work came from the
conjecture while it still remains unsettled. The author is grateful to Professor Zhan
for discussions and calling the author's attention to the paper \cite{Bo}.


\begin{thebibliography}{99}

\bibitem{AHJ}
T. Ando, R. Horn and C. Johnson,
The singular values of a Hadamard product: a basic inequality,
{\it Linear and Multilinear Algebra} {\bf 21} (1987), 345--365.

\bibitem{BS}
R.B. Bapat and V.S. Sunder,
On majorizaton and Schur products,
{\it Linear Algebra Appl.} {\bf 72} (1985), 107--117.

\bibitem{Bh}
R. Bhatia,
{\it Matrix Analysis}, Springer, New York, 1997.

\bibitem{Bo}
J.V. Bondar,
Comments on and complements to {\it Inequalities:\ Theory of Majorization and Its
Applications} by Albert W. Marshall and Ingram Olkin,
{\it Linear Algebra Appl.} {\bf 199} (1994), 115--130.

\bibitem{FH}
C.H. FitzGerald and R.A. Horn,
On fractional Hadamard powers of positive definite matrices,
{\it J. Math. Anal. Appl.} {\bf 61} (1977), 633--642.

\bibitem{Ha}
F. Hansen,
Functions of matrices with nonnegative entries,
{\it Linear Algebra Appl.} {\bf 166} (1992), 29--43.

\bibitem{Ho}
R.A. Horn,
The theory of infinitely divisible matrices and kernels,
{\it Trans. Amer. Math. Soc.} {\bf 136} (1969), 269--286.

\bibitem{HJ}
R.A. Horn and C.R. Johnson,
Hadamard and conventional submultiplicativity for unitarily invariant norms on
matrices,
{\it Linear and Multilinear Algebra} {\bf 20} (1987), 91--106.

\bibitem{MO}
A.W. Marshall and I. Olkin,
{\it Inequalities: Theory of Majorization and Its Applications},
Academic Press, New York, 1979.

\bibitem{O}
K. Okubo,
H\"older-type norm inequalities for Schur products of matrices,
{\it Linear Algebra Appl.} {\bf 91} (1987), 13--28.

\end{thebibliography}
\end{document}